\documentclass[reqno,10pt]{article}
\usepackage{hyperref}
\usepackage{amsmath}
\usepackage{amssymb}
\usepackage{mathrsfs}

\newtheorem{thm}{Theorem}[section]
\newtheorem{lemma}[thm]{Lemma}

\newtheorem{corol}[thm]{Corollary}
\newtheorem{propos}[thm]{Proposition}
\newtheorem{rema}{Remark}[section]
\def\bp{\begin{propos}}
\def\ep{\end{propos}}
\def\bt{\begin{thm}}
\def\et{\end{thm}}
\def\bco{\begin{corol}}
\def\eco{\end{corol}}
\def\bl{\begin{lemma}}
\def\el{\end{lemma}}
\def\br{\begin{rema}}
\def\er{\end{rema}}
\def\be{\begin{equation}}
\def\ee{\end{equation}}
\def\ba{\begin{array}}
\def\ea{\end{array}}
\def\bena{\begin{eqnarray}}
\def\eena{\end{eqnarray}}

\def\P{{\mathbb P}}
\def\E{{\mathbb E}}

\def\1{I}

\def\fF{{\mathscr F}}

\def\a{{\alpha}}
\def\D{{\Delta}}

\def\Var{\hb{Var}}
\def\b{\binom}

\def\QED{\hfill$\square$\vskip 3mm}

\def\Dp{\displaystyle}
\def\Df{\Dp\frac}
\def\hb{\hbox}
\def\({\left(}
\def\){\right)}

\begin{document}

\title{Phase Transition on the Degree Sequence of a Mixed Random Graph
Process\\[5mm]
 \footnotetext{Key
words and phrases: real-world networks; degree sequence; phase
transition; comparing argument.}}

\author{{Xian-Yuan Wu$^1$}\thanks{Supported in part
by the Natural Science Foundation of China},\ \ \ {Zhao
Dong$^2$}\thanks{Supported in part by the Natural Science Foundation
of China under grants 10671197 and 10721101},\ \ \ {Ke
Liu$^2$}\thanks{Supported in part by the Natural Science Foundation
of China under grants 60674082, 70221001 and 70731003.},\ \ \hb{ and
}{Kai-Yuan Cai}$^3$\thanks{Supported in part by the Natural Science
Foundation of China and MicroSoft Research Asia under grant
60633010}} \vskip 10mm
\date{}
\maketitle {\small \vskip-20mm \begin{center}
\begin{minipage}{14.5cm} \noindent\hskip -2mm$^1$School of
Mathematical Sciences, Capital Normal University,
Beijing, 100037, China. Email: \texttt{wuxy@mail.cnu.edu.cn}\\[-5mm]

\noindent\hskip -2mm$^{2}$Academy of Mathematics and System
Sciences, Chinese Academy of Sciences, Beijing, 100190, China.
Email: \texttt{dzhao@amss.ac.cn;}\ \texttt{kliu@amss.ac.cn}\\[-5mm]

\noindent\hskip -2mm$^3$Department of Automatic Control, Beijing
University of Aeronautics and Astronautics, Beijing, 100083, China.
Email: \texttt{kycai@buaa.edu.cn}
\end{minipage}
\end{center}
\vskip 2mm
\begin{center} \begin{minipage}{14.5cm}
{\bf Abstract}: This paper focuses on the problem of the degree
sequence for a mixed random graph process which continuously
combines the {\it classical} model and the BA model. Note that the
number of step added edges for the mixed model is random and
non-uniformly bounded. By developing a comparing argument, phase
transition on the degree distributions of the mixed model is
revealed: while the {\it pure} classical model possesses a {\it
exponential} degree sequence, the {\it pure} BA model and the mixed
model possess {\it power law} degree sequences. As an application of
the methodology, phase transition on the degree sequence of {\it
another} mixed model with {\it hard copying} is also studied,
especially, in the power law region, the inverse power can take any
value greater than 1.




\end{minipage}
\end{center}}

 \vskip 5mm
\section{Introduction and statement of the results}
\renewcommand{\theequation}{1.\arabic{equation}}
\setcounter{equation}{0}

Graph theory \cite{B,ER,G,H} is a rich research area that can be
traced back to the problem on the seven bridges in K\"onigsberg
considered by celebrated mathematician Euler in 1736. In 1950s,
Hungarian mathematicians Erd\"os and R\'enyi  extended the graph
theory into random environments and developed the classical theory
of random graphs. In their paper \cite{ER}, Erd\"os and R\'enyi
define the random graph $G_{n,M}$ (ER model) which consists of $n$
nodes and $M$ randomly chosen edges of the all $\Dp\b{n}{2}$
possible edges, and study the property of $G_{n,M}$ as $n\rightarrow
\infty$, with $M=M(n)$ a function of $n$. At the time when Erd\"os
and R\'enyi started their investigations of $G_{n,M}$, Gilbert
\cite{G} introduced a more fundamental random graph model $G_{n,p}$
as follows: Given $n$ nodes, each of the $\Dp\b n2$ distinct couples
of nodes is linked with an edge with probability $p$. For $M\sim
p\Dp\b n2$ as $n\rightarrow\infty$, the models $G_{n,M}$ and
$G_{n,p}$ are almost interchangeable and are subsequently called the
{\it classical} random graph models in the literature. Clearly, the
generation mechanism of the classical random graph is featured with
several characteristics. First, the number of nodes is given a {\it
priori} and keeps constant during the process of graph generation.
Second, the edges are generated in a random manner. Finally, each
edge is generated in an equal probability.

On the other hand, in recent years complex networks have drawn a lot
of attentions in disparate communities including statistical
mechanics, computer networks, control theory, among others
\cite{AB,ASBS,BKEMS,BKM,LLJ,N}. Various models involving random
factors have been proposed and investigated. Among them, the model
proposed by A.-L. Barab\'asi and R. Albert \cite{BA} (BA model) has
been well received and can be described as follows. A graph with
$n_0$ nodes and $m_0$ edges is given at the beginning. Then the
graph starts to evolve. At each time a new node with several new
edges is added to the graph. While all these new edges are linked
with the new node, the other node that links an edge of these new
edges is selected from the existing nodes according to the principle
of {\it preferential attachment}. Suppose that there are $n$ nodes
in the graph already, with $d_{x_i}$ being the degree of the $i^{\rm
th}$ node. The principle of preferential attachment asserts that the
$i^{\rm th}$ node is selected as the node that links one of the $m$
edges with probability ${d_{x_i}}/\sum_{i=1}^nd_{x_i}$. It is shown
that the degree distribution of the resulting graph obeys a power
law. Different from the generation mechanism of Erd\"os and R\'enyi,
for a random graph, the generation mechanism of BA is featured with
the following characteristics. First, the size of the graph in terms
of the number of nodes and edges is varying during the process of
graph generation. The graph tends to evolve. Second,
the added new edges are generated with unequal probabilities
according the principle of {\it preferential attachment}. Obviously,
the BA model can hardly be treated as an extension of the ER model.

A natural question is how to reconcile the ER theory of random
graphs and various models of complex networks and develop a coherent
or modern theory of random graph and complex networks, this forms
the {\it first motivation} of the present paper. As a useful step,
it should be interesting to combine the distinct features of the two
graph generation mechanisms described above and investigate various
properties of the resulting graph. In this paper we will first
introduce an evolving {\it classical} random graph model and then
modify this classical model according to the principle of {\it
preferential attachment}.

The ER model can be easily modified in an {\it evolving} way as
follows. Fix some constant $\mu>0$. Let's consider the following
process which generates a sequence of simple graphs
$\{G^0_t=(V_t,E_t)$, $t\geq 1\}$:

\vskip 2mm

{\it Time-Step 1.} {Let $G^0_1$ consists of vertices $x_0,x_1$ and
the edge $\langle x_0,x_1\rangle$. In general, $\langle u,v\rangle$
denotes the edge with endpoints $u,v$. }

\vskip 2mm {\it Time-Step $t\geq$ 2.} We add a vertex $x_t$ to
  $G^0_{t-1}$ and then add random edges incident with $x_t$: for any
  $0\leq i\leq t-1$, edge $\langle x_i,x_t\rangle$ is added
  independently with probability ${(\mu\wedge t)}/{t}$.
\vskip 2mm The process $\{G^0_t:t\geq 1\}$ defined above is called
{\it classical}, for edges are added in an equal probability at any
Time-Step, which coincides with the basic feature of ER model.

It is easily observed that the classical model $\{G^0_t:t\geq 1\}$
is not appropriate for studying real world networks also. Actually,
model $\{G^0_t:t\geq 1\}$ can be farther modified to the following
BA model $\{G_t=(V_t,E_t): t\geq 1\}$, which fits the first
motivation of us:

\vskip 2mm {\it Time-Step 1.} {Let $G_1$ consists of vertices
$x_0,x_1$ and the edge $\langle x_0,x_1\rangle$.}

\vskip 2mm {\it Time-Step $t\geq$ 2.} We add a vertex $x_t$ to
  $G_{t-1}$ and then add random edges incident with $x_t$: for any
  $0\leq i\leq t-1$, edge $\langle x_i,x_t\rangle$ is added
  independently with probability $\frac{\mu d_{x_i}(t-1)}{2e_{t-1}}\wedge 1$, where
  $d_{x_i}(t-1)$ be the degree of $x_i$ in $G_{t-1}$ and $e_{t-1}=|E_{t-1}|$.

\vskip 2mm The {\it second motivation} for us to consider the above
process $\{G_t: t\geq 0\}$ is to model the www-{\it typed}
real-world networks properly. We say a real-world network is of
www-{\it typed}, if the following holds
\begin{enumerate}
  \item Excepting for all the isolated vertices (nodes), the network
  has only one connected component;
  \item There is no loop and multi-edge in the network;
  \item While a new vertex (node) is added, the number of added new edges (links) between it and the existing vertices is finite but unbounded; and
  \item Edges (links) are added in the {\it preferential
  attachment} manner.
\end{enumerate}

Actually, to model the real world networks by random complex graphs,
many new models (deferring from the ER model) have already been
introduced. By studying complex graphs, various topological
properties such as degree-distribution \cite{BA,BO,BRST,FFF},
diameter \cite{ABJ,ASBS,BR2,SAB}, clustering \cite{BR,N}, stability
\cite{B,B2,BR3} and spectral gap \cite{ACL} of these real-world
networks have been presented. One of the most basic properties of
real-world networks is the power law degree distribution, many new
models have been introduced to explain the underlying causes for the
emergence of power law degree distributions. This can be observed in
the `LCD model' \cite{BR2}; the generalization of `LCD model' due to
Buckley and Osthus \cite{BO}; `copying' models of Kumar {\it et al}.
\cite{KRRS}; `hard copying' models of Wu {\it et al.} \cite{NWC};
the general models defined by Copper and Frieze \cite{CF}; the
growth-deletion models of Copper, Frieze and Vera \cite{CFV}, Chung
and Lu \cite{CL} and Wu {\it et al}. \cite{Wu} {\it etc}. The main
difference between our model and those introduced in
\cite{BO,BR2,CF,CFV,CL,KRRS} and \cite{Wu} is that, in our setting,
the number of step added edges is random and non-uniformly bounded.
Note that the `hard copying' model introduced in \cite{NWC} is also
a model with non-uniformly bounded edge addition. Obviously, the
model $\{G_t:t\geq 1\}$ seems to be a more proper candidate for
modeling the www-{\it typed} real-world networks.

Now, Let $D_k(t)$ be the number of vertices with degree $k\geq 0$ in
$G_t$ and let $\overline{D}_k(t)$ be the expectation of $D_k(t)$.
Note that, in this paper, for any kind of random graph process, we
always denote $D_k(t)$ the number of vertices with degree $k\geq 0$
and $\overline{D}_k(t)$ its expectation.

The first result of this paper is about BA model $\{G_t=(V_t,E_t):
t\geq 1\}$, it follows as

\bt\label{th1} For any $0<\mu\leq 2$, there exists positive
constants $C_1$ and $C_2$ such that
 \be\label{1.1}
C_1k^{-3}\leq\liminf_{t\rightarrow\infty }\Df{\overline
{D}_k(t)}{t}\leq \limsup_{t\rightarrow\infty }\Df{\overline
{D}_k(t)}{t}\leq C_2 k^{-3}\ee for all $k\geq 1$.

\et

\br In this paper, the condition $0<\mu\leq 2$ is purely technical,
and it is conjectured that our results hold for any $\mu>0$.
\er

By definition, excepting for the isolated vertices, $G_t$ contains a
unique connected component, we call it the {\it giant component} of
$G_t$. Denote by $C_t$ the giant component. The following is our
result on $\E(|C_t|)$, the mean size of $C_t$.

\bt\label{th2} Assume that $0<\mu\leq 2$. Then for any small enough
$\nu>0$, we have
 \be\label{1.2}
\E(|C_t|)=(1-e^{-\mu})t+O(t^{\frac 1{2-\nu}}).\ee Note that the
hidden constant in $O(t^{\frac 1{2-\nu}})$ only depends on $\nu$.
\et

Now, we present a mixed model which continuously combines the {\it
classical} model $\{G^0_t:t\geq 1\}$ and the above BA model
$\{G_t:t\geq 1\}$. Fix some constants $0\leq \a\leq 1$ and $\mu,\
\zeta>0$. Define random graph process $\{G^\a_t=(V_t,E_t): t\geq
1\}$ as follows.

\vskip 2mm

{\it Time-Step 1.} {Let $G^\a_1$ consists of vertices $x_0,x_1$ and
the edge $\langle x_0,x_1\rangle$.}

\vskip 2mm {\it Time-Step $t\geq$ 2.} We add a new vertex $x_t$ to
$G^\a_{t-1}$ and then \vskip 3mm
\begin{enumerate}
  \item  with probability $\a$, we add random edges incident with $x_t$ in the {\it preferential
attachment} manner: for any $0\leq i\leq t-1$, edge $\langle
x_i,x_t\rangle$ is added independently with probability $\frac{\mu
d_{x_i}^\a(t-1)}{2e_{t-1}}\wedge1$, where $d_{x_i}^\a(t-1)$ be the
degree of $x_i$ in $G^\a_{t-1}$;
  \item with probability $1-\a$, we add random edges incident with $x_t$ in the {\it classical
} manner: for any $0\leq i\leq t-1$, edge $\langle x_i,x_t\rangle$
is added independently with probability $(\zeta\wedge t) /t$.
\end{enumerate}

It is straightforward to generalize the approach developed for
Theorem~\ref{th1} to prove the following corollary for $\{G^\a_t:
t\geq 1\}$, $0<\a<1$:

\bco\label{th3} For any $0<\a<1$, $0<\mu\leq 2$ and $\zeta>0$, there
exists positive constants $C^\a_1$ and $C^\a_2$ such that
 \be\label{1.3}
C^a_1k^{-\beta}\leq\liminf_{t\rightarrow\infty }\Df{\overline
{D}_k(t)}{t}\leq \limsup_{t\rightarrow\infty }\Df{\overline
{D}_k(t)}{t}\leq C^\a_2 k^{-\beta}\ee for all $k\geq 1$, where
$\beta=1+2\(1+\Dp\frac{(1-\a)\zeta}{\a\mu}\)$.\eco

\br\label{r4} At any Time-Step $t>\zeta$, the mean number of added
new edges is $\xi:=\a\mu+(1-\a)\zeta$ and
$\frac{(1-\a)\zeta}{\a\mu}$ be the limit ratio of the number of the
two kinds of edges in $G^\a_t$.\er

In the case of $\a=0$, we get the {\it classical} process
$\{G^0_t:t\geq 1\}$ parameterized by $\zeta>0$. Just as one expects,
the model $\{G^0_t:t\geq 1\}$ possesses a {\it classical
(exponential) } degree sequence as

\bco\label{th4} For random graph process $\{G^0_t:t\geq 1\}$, there
exists positive constants $C^0_1$ and $C^0_2$ such that
 \be\label{1.4}
C^0_1\(\frac{\zeta}{1+\zeta}\)^k\leq\liminf_{t\rightarrow\infty
}\Df{\overline {D}_k(t)}{t}\leq \limsup_{t\rightarrow\infty
}\Df{\overline {D}_k(t)}{t}\leq C^0_2 \(\frac{\zeta}{1+\zeta}\)^k
\ee for all $k\geq 0$.\eco

Theorems~\ref{th1} and Corollaries \ref{th3} and \ref{th4} exhibit a
phase transition on the degree distributions of the mixed model
$\{G^\a_t:t\geq 1\}$ while $\a$ varies from $0$ to $1$. Note that
phase transition on degree distributions of random graph process is
first studied in the recent work \cite{Wu} of Wu {\it et al.}. More
precisely, \cite{Wu} introduced a model with edge deletions and
showed that, while a relevant parameter varies, the model exhibits
{\it power law} degree distribution, a special degree distribution
lying {\it between power law and exponential}, and {\it exponential}
degree distribution in turn. A numerical investigation to phase
transition on degree distributions of networks can be founded in
reference \cite{ZJW}.

The rest of the paper is organized as follows. In Section 2, we give
some useful estimates to $e_t$, the number of edges in $G_t$. In
section 3, we bound the maximum degree of vertex in $G_t$, and then
prove Theorem~\ref{th2}. In Section 4, we establish the recurrence
for $\overline{D}_k(t)$, then solve the recurrence by using a
compare argument, and finally finish the proof of Theorem~\ref{th1}.
In Section 5, we adopt the comparing argument developed in Section 4
to prove Corollaries~\ref{th3} and \ref{th4}. In Section 6, we apply
the comparing argument to study the phase transition on the degree
sequence of a mixed model with {\it hard copying}.

\section{Estimates for $e_t$}
\renewcommand{\theequation}{2.\arabic{equation}}
\setcounter{equation}{0} In this section we give some lemmas for
$e_t$, which will play important roles in the proofs of our main
results.

We first consider the increments of $e_t$. Let $a_t=e_{t+1}-e_t$ and
$\{\fF_t:t\geq 1\}$ be the natural $\sigma$-flow generated by
process $\{G_t:t\geq 1\}$. Then \bl\label{l2.1}For all $t\geq 1$, we
have \be\label{2.0}\E(a_t\mid\fF_t)=\mu\ee and
\be\label{2.1}\E(a_t^k\mid\fF_t)\leq (\mu\vee 1)^kk!\ee for $k\geq
2$.\el

{\it Proof}: Let $\{p_i:1\leq i\leq n\}$, $n\geq 2$, be a serial of
positive numbers satisfying $p_i\leq \frac 12$, $\sum_{i=1}^np_i=1$,
and let $\{X_i,1\leq i\leq n\}$ be the independent random variables
with
$$\P(X_i=1)=\mu p_i=1-\P(X_i=0).$$ Let $X=\Dp\sum_{i=1}^nX_i$. Clearly, to prove the
lemma, it suffices to prove that $$\E(X)=\mu\ \ \hb{and} \ \
\E(X^k)\leq (\mu\vee 1)^k\times k!\ \ \forall\ \ k\geq 2.$$

For $k=1$, it is straightforward to see that $\E(X^k)=\E(X)=\mu\leq
\mu\vee 1$. Assume that $\E(X^{m})\leq (\mu\vee 1)^{m}\times m!$ for
some $m\geq 1$, then \bena
&&\E(X^{m+1})=\E(\sum_{i=1}^nX_i)^{m+1}=\sum_{i_1=1}^n\cdots\sum_{i_m=1}^n\sum_{i_{m+1}=1}^n\E(X_{i_1}\ldots
X_{i_m}X_{i_{m+1}})\nonumber\\[3mm]
&&= \sum_{i_1=1}^n\cdots\sum_{i_m=1}^n\(\sum_{i_{m+1}\in\{i_1,\ldots, i_m\}}\E(X_{i_1}\ldots
X_{i_m})+\sum_{i_{m+1}\notin\{i_1,\ldots, i_m\}}\E(X_{i_1}\ldots
X_{i_m})\E(X_{i_{m+1}})\)\nonumber\\[3mm]
&&\leq \sum_{i_1=1}^n\cdots\sum_{i_m=1}^n\(m\E(X_{i_1}\ldots
X_{i_m})+\mu\E(X_{i_1}\ldots X_{i_m})\)\nonumber\\[3mm]&& \leq (m+1)(\mu\vee
1)\sum_{i_1=1}^n\cdots\sum_{i_m=1}^n\E(X_{i_1}\ldots
X_{i_m})\nonumber\\[3mm]&&=(m+1)(\mu\vee 1)\E(X^m)\leq (\mu\vee
1)^{m+1}\times(m+1)!.\nonumber\eena Thus we finish the proof by
induction.\QED

Now, define $Y_t=e_t-\mu t$ for $t\geq 1$, then, by the definition
of $G_t$, $\{Y_t:t\geq 1\}$ forms a martingale with respect to
$\{\fF_t:t\geq 1\}$.

\bl\label{l2.2}There exists some constant $c_1>0$ such
that\be\label{2.2} \P(|e_t-\mu t|\geq t^{4/5})\leq c_1 t^{-3/5}\ee
for all $t\geq 1$.\el

{\it Proof}: By the property of martingale, first, we have
\be\label{2.3}\E(Y_t-Y_1)^2=\E\(\sum_{i=1}^{t-1}(Y_{i+1}-Y_i)\)^2=\sum_{i=1}^{t-1}\E(Y_{i+1}-Y_i)^2=\sum_{i=1}^{t-1}\Var(a_i).\ee
Then, by Lemma~\ref{l2.1}
\be\label{2.4}\E(Y_t-Y_1)^2=\sum_{i=1}^{t-1}\Var(a_i)\leq
\(2(\mu\vee 1)^2-\mu^2\)(t-1).\ee Finally, using the relation that
$\E(Y_t^2)=\E(Y_t-Y_1)^2+(1-\mu)^2$ and the Markov's inequality, we
have \[\P(|e_t-\mu t|\geq t^{4/5})\leq\Dp\frac
{\E(Y^2_t)}{t^{8/5}}\leq \Df{(2(\mu\vee
1)^2-\mu^2)(t-1)+(1-\mu)^2}{t^{8/5}}\leq c_1 t^{-3/5}\] for some
constant $c_1>0$.\QED

\bl\label{l2.3} For any $\nu>0$, there exists constants $c_2, c_3>0$
such that\be\label{2.5}\P(|e_t-\mu t|\geq \nu t)\leq c_2e^{-c_3t}\ee
for all $t\geq 1$.\el

{\it Proof}: By Lemma~\ref{l2.1}, for small $\lambda>0$, we have
\[\E(e^{\lambda a_t}\mid \fF_t)=1+\lambda\mu+O(\lambda^2),\]then
\[\E(e^{\lambda e_{t+1}})=\E\(\E(e^{\lambda e_t+\lambda
a_t}\mid\fF_t)\)=\E\(e^{\lambda e_t}\E(e^{\lambda
a_t}\mid\fF_t)\)=(1+\lambda\mu+O(\lambda^2))\E(e^{\lambda e_t}).\]
This implies that
\[\E(e^{\lambda e_t})=(1+\lambda\mu+O(\lambda^2))^{t-1}\E(e^{\lambda e_1})=\frac{e^\lambda}{1+\lambda\mu+O(\lambda^2)}
\exp\{\ln(1+\lambda\mu+O(\lambda^2))t\}.\] For given $\nu>0$, take
$\lambda>0$ small enough such that
\[c'_3:=(\mu+\nu)\lambda-\ln(1+\lambda\mu+O(\lambda^2))>0.\] Taking
$c'_2={e^\lambda}/{\(1+\lambda\mu+O(\lambda^2)\)}$, we have
\be\label{2.7}\P(e_t\geq (\mu+\nu)t)\leq \E(e^{\lambda
e_t})e^{-(\mu+\nu)\lambda t}\leq c'_2e^{-c'_3t}.\ee Similarly, for
some $c''_2,c''_3>0$, we have \be\label{2.8}\P(e_t\leq
(\mu-\nu)t)\leq e^{(\mu-\nu)\lambda t}\E(e^{-\lambda e_t})\leq
c''_2e^{-c''_3t}.\ee The lemma follows from (\ref{2.7}) and
(\ref{2.8}).\QED

\section{Bounding the degree and the proof of Theorem~\ref{th2}}

\renewcommand{\theequation}{3.\arabic{equation}}
\setcounter{equation}{0}

For times $s$ and $t$ with $0\leq s\leq t$, $t\geq 1$, let
$d_{x_s}(t)$ be the degree of vertex $x_s$ in $G_t$. In this
section, we will concentrate on the upper bound of $d_{x_s}(t)$ and
then prove Theorem~\ref{th2}.

We say an event happens {\it quite surely} (qs) if the probability
of the complimentary set of the event is $O(t^{-K})$ for any $K>0$.

The following is our bounding for $d_{x_s}(t)$. As noted in
\cite{Wu}, our result will depend on Lemma~\ref{l2.3}, the
exponential inequality for $e_t$.

\bl\label{l3.1} For small $\nu>0$ and $1\leq s\leq t$, we have
\be\label{3.1}d_{x_s}(t)\leq(t/s)^{\frac1{2-\nu}}(\log t)^3 \ \ \
\mbox{qs}.\ee\el

{\it Proof}: Let $X_\tau=d_{x_s}(\tau)$ for $\tau=s,s+1,\ldots,t$.
Conditional on $X_\tau=x$ and $e_\tau$, we have
\be\label{3.2}X_{\tau+1}= x+B\(1,\Dp\frac{\mu x}{2e_\tau}\),\ee
where $B\(1,\Dp\frac{\mu x}{2e_\tau}\)$ be the $\{0,1\}$-valued
random variable with $\P\(B\(1,\Dp\frac{\mu
x}{2e_\tau}\)=1\)=\Dp\frac{\mu x}{2e_\tau}.$

Lemma~\ref{l3.1} follows immediately from (\ref{2.5}), (\ref{3.2})
and a standard argument which can be found in the proof of Lemma 2.1
in \cite{Wu}. \QED

\br\label{r1}Because $d_{x_0}(t)$ and $d_{x_1}(t)$ are same
distributed, Lemma~\ref{l3.1} implies that \[d_{x_0}(t)\leq t^{\frac
1{2-\nu}}(\log t)^3,\ \ \ qs.\] \er

Now, based on Lemma~\ref{l3.1}, we prove Theorem~\ref{th2} as
follows.

{\it Proof of Theorem~\ref{th2}}: To prove Theorem~\ref{th2}, it
suffices to show that \be\label{3.9}\E(|V_t\setminus
C_t|)=e^{-\mu}t+O(t^{\frac 1{2-\nu}}).\ee Denote by $\Delta_t$ the
maximal degree in $G_t$. By Lemma~\ref{l2.3}, Lemma~\ref{l3.1} and
Remark~\ref{r1}, we have \be\label{3.10}\Df{\Delta_t}{e_t}\leq L
t^{\frac 1{2-\nu}-1},\ \ \ qs \ee where $L$ be a constant
independent of $t$.

For large $t$, let's consider the probability $\P(a_t=0)$, recall
that $a_t=e_{t+1}-e_t$ be the increment of $e_t$ at Time-Step $t+1$.
By equation (\ref{3.10}), we have
\bena&&\P(a_t=0)=\E(\1_{a_t=0})=\E(\E(\1_{a_t=0}\mid\fF_t))\nonumber\\[2mm]
&&\hskip-4mm=\E\(\E(\1_{a_t=0}\mid\fF_t)\left|
\Df{\Delta_t}{e_t}\leq L t^{\frac
1{2-\nu}-1}\right.\)\P\(\Df{\Delta_t}{e_t}\leq L t^{\frac
1{2-\nu}-1}\)\nonumber\\[3mm]&& +\E\(\E(\1_{a_t=0}\mid\fF_t)\left|
\Df{\Delta_t}{e_t}> L t^{\frac
1{2-\nu}-1}\right.\)\P\(\Df{\Delta_t}{e_t}> L t^{\frac 1{2-\nu}-1}\)
\nonumber\\[3mm]&&
\hskip-4mm=\E\(\E(\1_{a_t=0}\mid\fF_t)\left| \Df{\Delta_t}{e_t}\leq
L t^{\frac 1{2-\nu}-1}\right.\)+O(t^{-10}).\label{3.11}\eena

The term $\E(\1_{a_t=0}\mid\fF_t)$ can be expressed as
\bena&&\E(\1_{a_t=0}\mid\fF_t)=\prod_{i=0}^t(1-\frac{\mu
d_{x_i}(t)}{2e_t})=\exp\left\{\sum_{i=0}^t\log\(1-\frac{\mu
d_{x_i}(t)}{2e_t}\)\right\}\nonumber\\[3mm]&&
\hskip-3mm=\exp\left\{-\sum_{i=0}^t\frac{\mu
d_{x_i}(t)}{2e_t}+O\(\sum_{i=0}^t\(\frac{\mu
d_{x_i}(t)}{2e_t}\)^2\)\right\}
=e^{-\mu}+O\(\frac{\Delta_t}{e_t}\),\label{3.11'}\eena hence,
\be\label{3.12}\E\(\E(\1_{a_t=0}\mid\fF_t)\left|
\Df{\Delta_t}{e_t}\leq L t^{\frac
1{2-\nu}-1}\right.\)=e^{-\mu}+O\(t^{\frac 1{2-\nu}-1}\).\ee Thus,
(\ref{3.11}) and (\ref{3.12}) imply that
\be\label{3.13}\P(a_t=0)=e^{-\mu}+O\(t^{\frac 1{2-\nu}-1}\).\ee

Now, by the definition of $G_t$, we
have\be\label{3.14}\E(|V_t\setminus
C_t|)=\sum_{s=2}^t\P(d_{x_s}(t)=0)=\sum_{s=1}^{t-1}\P(a_{s}=0),\ee
equation (\ref{3.9}) follows immediately from (\ref{3.13}) and
(\ref{3.14}).\QED

\br\label{r2}For any $t\geq 1$, we have
\bena&&\P(a_t=0)=\E(\E(\1_{a_t=0}\mid\fF_t))=\E\(\prod_{s=0}^t\(1-\frac{\mu
d_{x_s}(t)}{2e_t}\)\)\nonumber\\&&\leq
\E\(\prod_{s=0}^t\exp\left\{-\frac{\mu
d_{x_s}(t)}{2e_t}\right\}\)=e^{-\mu}.\label{3.16}\eena Furthermore,
equation (\ref{3.13}) implies that
$\Dp\lim_{t\rightarrow\infty}\P(a_t=0)=e^{-\mu}$.

For the probability $\P(a_t=1)$, using (\ref{3.10}) again, the same
arguments as in (\ref{3.11}-\ref{3.12}) imply that
\be\label{3.15}\lim_{t\rightarrow\infty}\P(a_t=1)=\mu e^{-\mu}.\ee
\er

\section{The comparing Approach and The proof of Theorem~\ref{th1}}

\renewcommand{\theequation}{4.\arabic{equation}}
\setcounter{equation}{0}

In this Section, we develop a comparing approach to prove
Theorem~\ref{th1}. We first follow the basic procedures in
\cite{CFV} to establish the recurrence for $\overline D_k(t)$. By
the definition of $G_t$, first of all, we have $D_0(1)=0$,
$D_1(1)=2$ and $D_k(t)=0$ for all $k,t$ with $k>t\geq 1$.

Now, put $D_{-1}(t)=0$ for all $t\geq 1$. For $t+1\geq k\geq 0$ and
$t\geq 1$, we have
\be\label{4.1}\E(D_k(t+1)\mid\fF_t)=D_k(t)+\left(-\frac{k\mu
D_k(t)}{2e_t}+\frac{(k-1)\mu
D_{k-1}(t)}{2e_t}\right)+\E(\1_{a_t=k}\mid\fF_t).\ee Taking
expectation and then using the basic inequality
$$e_t\leq\sum_{s=1}^ts=\Df{t(t+1)}{2}$$ and the estimations given in Lemmas~\ref{l2.2} and~\ref{l2.3},
(\ref{4.1}) implies that\be\label{4.1'}\overline D_k(t+1)=\overline
D_k(t)+\Dp\frac{k-1}{2}\frac {\overline
D_{k-1}(t)}{t}-\Dp\frac{k}{2}\frac {\overline
D_{k}(t)}{t}+O(t^{-1/5})+f_k(t),\ee where $f_k(t)=\P(a_t=k)$. Note
that term $O(t^{t^{-1/5}})$ is independent of $k$. We get the
recurrence for $\overline D_k(t)$ as:
\be\label{4.2}\left\{\begin{array}{rl}&\hskip-4mm\overline
D_k(t+1)=\overline D_k(t)+\Dp\frac{k-1}{2}\frac {\overline
D_{k-1}(t)}{t}-\Dp\frac{k}{2}\frac
{\overline D_{k}(t)}{t}+O(t^{-1/5})+f_k(t),\\[4mm]&\hskip 80mm  t+1\geq k\geq 0,\ t\geq 1;\\[4mm]
&\hskip-4mm\overline D_0(1)=0;\hskip 3mm\overline D_1(1)=2;\hskip
3mm\overline D_k(t)=0,\ k> t\geq 1;\hskip 3mm\overline D_{-1}(t)=0,\
t\geq 1.\end{array}\right.\ee

To solve the recurrence (\ref{4.2}), we need a comparing argument.
Note that the recurrence as (\ref{4.2}) with $\{f_k(t)\}$ replaced
by a serial of constants can be solved directly by the method
developed in \cite{CF}, \cite{CFV} and \cite{Wu}. Let
\[F_k(t):=\overline D_k(t+1)-\overline D_k(t)-\Dp\frac{k-1}{2}\frac
{\overline D_{k-1}(t)}{t}+\Dp\frac{k}{2}\frac {\overline
D_{k}(t)}{t}-f_k(t).\] Obviously, $F_k(t)$ is a determined (or
known!) function in $k$ and $t$ satisfying
\be\label{4.3}|F_k(t)|\leq Rt^{-1/5},\ \ \forall\ k\geq 0,\ t\geq
1.\ee For $k\geq 0$, define
\[A_k(t)=\left\{\ba{lll}&\hskip-4mm F_k(t),&\hb{ if  }t\geq k,\\[3mm]
&\hskip-4mm F_k(t)+f_k(t),&\hb{ if  }t\leq k-1;\ea\right.\hb{  and }\hskip 5mm g_k(t)=\left\{\ba{lll}&\hskip-4mm f_k(t),&\hb{ if  }t\geq k,\\[3mm]
&\hskip-4mm 0,&\hb{ if  }t\leq k-1.\ea\right.\] Then, (\ref{4.2})
can be rewritten as
\be\label{4.6}\left\{\begin{array}{rl}&\hskip-4mm\overline
D_k(t+1)=\overline D_k(t)+\Dp\frac{k-1}{2}\frac {\overline
D_{k-1}(t)}{t}-\Dp\frac{k}{2}\frac
{\overline D_{k}(t)}{t}+A_k(t)+g_k(t),\\[4mm]&\hskip 80mm  t+1\geq k\geq 0, \ t\geq 1;\\[4mm]
&\hskip-4mm\overline D_0(1)=0;\hskip 3mm\overline D_1(1)=2;\hskip
3mm\overline D_k(t)=0,\ k> t\geq 1;\hskip 3mm\overline D_{-1}(t)=0,\
t\geq 1.\end{array}\right.\ee By the fact that $f_k(t)=0$ for $t\leq
k-2$ and $f_k(k-1)=\P(a_{k-1}=k)\leq \mu k^{-1}$ for $k\geq 2$,
similar to (\ref{4.3}), we have for some $ R_1>0$ \be\label{4.6'}
|A_k(t)|\leq R_1t^{-1/5},\ \ \forall\ k\geq 0,\ t\geq 1.\ee

In the rest of this section, we will try to solve the recurrence
(\ref{4.6}) for any given function serial $\{A_k(t)\}$ satisfying
(\ref{4.6'}). The lack of the existence of such limit as
$\lim_{t\rightarrow\infty} f_k(t)$ makes it difficult to solve
(\ref{4.6}) directly. In fact, to solve (\ref{4.6}) by the known
argument developed in \cite{CF}, \cite{CFV} and \cite{Wu}, we not
only need the existence of such limits, but also need a uniform
speed faster than $t^{-\epsilon}$, $\epsilon>0$, of the
corresponding convergence. But this seems impossible (see the proof
of Corollary~\ref{th4}), we have to develop a new method to study
$\overline D_k(t)$.

By Remark~\ref{r2}, $\lim_{t\rightarrow\infty}\P(a_t=0)=e^{-\mu}$,
then, for some constant $\rho>0$, \be\label{4.0}\P(a_t=0)\geq
\rho>0,\ \ \forall \ t\geq 1.\ee
For $k\geq 0$, let \be\label{4.29}\psi_k=\left\{\ba{rll}&\hskip-4mm0,\hskip3mm &k\geq 1,\\[5mm]
&\hskip-4mm\rho,\hskip3mm &k=0;\ea\right.\hb{ and}\hskip3mm
\varphi_k=\left\{\ba{rll}&\hskip-4mm Ck^{-4},\hskip3mm
&k\geq 1,\\[5mm]
&\hskip-4mm e^{-\mu},\hskip3mm &k=0,\ea\right.\ee with $C=(\mu\vee
1)^4\times 4!$. Define
\be\label{4.30}\psi_k(t)=\left\{\ba{rll}&\hskip-4mm0,\hskip3mm
&k\geq 1,\ t\geq
1,\\[5mm]
&\hskip-4mm\psi_k,\hskip3mm &k=0,\ t\geq 1;\ea\right.\hb{
and}\hskip3mm
\varphi_k(t)=\left\{\ba{rll}&\hskip-4mm\varphi_k,\hskip3mm
&t\geq k,\\[5mm]
&\hskip-4mm0,\hskip3mm &1\leq t<k.\ea\right.\ee By Lemma~\ref{l2.1},
equation (\ref{3.16}) and the Markov's inequality, we have
\be\label{4.7}\psi_k(t)\leq g_k(t)\leq \varphi_k(t),\ \ \forall\
k\geq 0,\ t\geq 1.\ee Now, with $g_k(t)$ in (\ref{4.6}) replaced by
$\psi_k(t)$ and $\varphi_k(t)$ respectively, we get the following
recurrences for $\tilde D_k(t)$ and $\hat D_k(t)$:

\be\label{4.8}\left\{\begin{array}{rl}&\hskip-4mm\tilde
D_k(t+1)=\tilde D_k(t)+\Dp\frac{k-1}{2}\frac {\tilde
D_{k-1}(t)}{t}-\Dp\frac{k}{2}\frac
{\tilde D_{k}(t)}{t}+A_k(t)+\psi_k(t),\hskip 2mm t+1\geq k\geq 0,\ t\geq 1;\\[4mm]
&\hskip-4mm\tilde D_0(1)=0;\hskip 3mm\tilde D_1(1)=2;\hskip
3mm\tilde D_k(t)=0,\ k> t\geq 1;\hskip 3mm\tilde D_{-1}(t)=0,\ t\geq
1;\end{array}\right.\ee

\vskip2mm \be\label{4.9}\left\{\begin{array}{rl}&\hskip-4mm\hat
D_k(t+1)=\hat D_k(t)+\Dp\frac{k-1}{2}\frac {\hat
D_{k-1}(t)}{t}-\Dp\frac{k}{2}\frac
{\hat D_{k}(t)}{t}+A_k(t)+\varphi_k(t),\hskip 2mm t+1\geq k\geq 0,\ t\geq 1;\\[4mm]
&\hskip-4mm\hat D_0(1)=0;\hskip 3mm\hat D_1(1)=2;\hskip 3mm\hat
D_k(t)=0,\ k> t\geq 1;\hskip 3mm\hat D_{-1}(t)=0,\ t\geq
1.\end{array}\right.\ee

We first give the following comparing lemma to show that $\tilde
D_k(t)$ and $\hat D_k(t)$ are lower and upper bounds for $\overline
D_k(t)$ respectively.

\bl\label{l4.1}{ \rm \bf [Comparing Lemma]} Assume that $\tilde
D_k(t)$ and $\hat D_k(t)$ be the solutions of (\ref{4.8}) and
(\ref{4.9}) respectively. Then \be\label{4.10}\tilde D_k(t)\leq
\overline D_k(t)\leq \hat D_k(t),\ \ \forall\ k\geq-1,\ t\geq 1.\ee
\el

{\it Proof}: We only prove the first inequality in (\ref{4.10}), the
situation for the second one is the same. Firstly, noticing that
$\tilde D_{-1}(t)=\overline D_{-1}(t)=0$ for all $t\geq 1$, we have
\[\tilde D_0(t+1)=\tilde D_0(t)+A_0(t)+\psi_0(t)\]
and
\[\overline D_0(t+1)=\overline D_0(t)+A_0(t)+g_0(t)\]
for all $t\geq 1$. This, together with the fact that $\tilde
D_0(1)=\overline D_0(1)=0$ and the inequality (\ref{4.7}), implies
\be\label{4.11}\tilde D_0(t)\leq\overline D_0(t),\ \ \forall\ t\geq
1.\ee Secondly, by the fact that $\tilde D_{k+1}(k)=\overline
D_{k+1}(k)=\psi_{k+1}(k)=g_{k+1}(k)=0$ for all $k\geq 1$, we have
\be\label{4.12'}\tilde D_{k+1}(k+1)=\frac 12\tilde
D_{k}(k)+A_{k+1}(k)\ee and
\[\overline D_{k+1}(k+1)=\frac 12\overline D_{k}(k)+A_{k+1}(k)\]
for all $k\geq 1$. This, together with the initial condition $\tilde
D_1(1)=\overline D_1(1)=2$, implies that \be\label{4.12}\tilde
D_k(k)=\overline D_k(k),\ \ \forall\ k\geq 1.\ee

Suppose we have proved that for some $m\geq 0$,
\be\label{4.13}\tilde D_k(k+m)\leq \overline D_k(k+m),\ \ \forall\
k\geq 1.\ee If we can prove \be\label{4.14}\tilde D_k(k+(m+1))\leq
\overline D_k(k+(m+1)),\ \ \forall\ k\geq 1,\ee then we get the
lemma by induction.

By (\ref{4.7}) and (\ref{4.13}), (\ref{4.14}) can be easily proved
by induction. The details are omitted. \QED

Now we begin to solve (\ref{4.8}) and (\ref{4.9}). We introduce two
recurrences with respect to $\{\psi_k\}$ and $\{\varphi_k\}$ as
follows: \be\label{4.15}\left\{\ba{rl}&\hskip-4mm\tilde
d_k=\Dp\frac{k-1}{2}\tilde
d_{k-1}-\Dp\frac k2\tilde d_k+\psi_k,\ \ \ k\geq 0,\\[5 mm]
&\hskip-4mm\tilde d_{-1}=0;\ea\right.\ee
\be\label{4.16}\left\{\ba{rl}&\hskip-4mm\hat
d_k=\Dp\frac{k-1}{2}\hat
d_{k-1}-\Dp\frac k2\hat d_k+\varphi_k,\ \ \ k\geq 0,\\[5mm]
&\hskip-4mm\hat d_{-1}=0.\ea\right.\ee The following Lemma show that
(\ref{4.15}) and (\ref{4.16}) are good approximation to (\ref{4.8})
and (\ref{4.9}) respectively.

\bl\label{l4.2}Assume that $\{\tilde D_k(t): k\geq -1,t\geq 1\}$
(resp. $\{\hat D_k(t): k\geq -1,t\geq 1\}$)
 be the solution of recurrence (\ref{4.8}) (resp. (\ref{4.9}))and $\{\tilde d_k:k\geq -1\}$ (resp. $\{\hat d_k:k\geq -1\}$) be the solution
of (\ref{4.15})(resp. (\ref{4.16})). If $\tilde d_k\leq C/{k}$
(resp. $\hat d_k\leq C/{k}$) for $k>0$ and some constant $C$, then
there exists constant $M_1$ (resp. $M_2$) such that\be\label{4.17}
\left|\tilde D_k(t)-t\tilde d_k\right|\leq M_1t^{4/5}\ \ (resp.
\left|\hat D_k(t)-t\hat d_k\right|\leq M_2t^{4/5}) \ee for all
$k\geq -1$ and $t\geq 1$.\el

{\it Proof of Lemma~\ref{l4.2}}: By using the fact that $\tilde
D_k(t)=0$ (resp. $\hat D_k(t)=0$) for $k>t\geq 1$ and the condition
$\tilde d_k\leq C/{k}$ (resp. $\hat d_k\leq C/{k}$), it is
straightforward to prove Lemma~\ref{l4.2} by induction (in $t$).
Note that our inductive hypothesis is \[\tilde{\cal H}_t:\ \
|\tilde\Theta_k(t)|\leq M_1t^{4/5}\ \ \hb{for all}\ \ k\ge -1. \ \
(resp.\ \hat{\cal H}_t:\ \ |\hat\Theta_k(t)|\leq M_2t^{4/5}\ \
\hb{for all}\ \ k\ge -1.)\] For details, one may refer to \cite{Wu}
(the proof of Lemma 2.2).\QED

Now, we finish the proof of Theorem~\ref{th1} as follows.

{\it Proof of Theorem~\ref{th1}}: For any given constant number
serial $\{\phi_k: k\geq 0\}$, the recurrence in $k$ with the form
\[\left\{\ba{rl}&\hskip-4mm d_k=\Dp\frac{k-1}{2}
d_{k-1}-\Dp\frac k2 d_k+\phi_k,\ \ \ k\geq 0,\\[5 mm]
&\hskip-4mm d_{-1}=0,\ea\right.\] can be directly solved as:
$d_{-1}=0$, $d_0=\phi_0$, $d_1=\frac 23 \phi_1$ and \be\label{4.27}
d_k=\sum_{j=1}^{k}\frac {2j(j+1)}{k(k+1)(k+2)}\phi_j=\frac
1{k(k+1)(k+2)}\sum_{j=1}^{k} {2j(j+1)}\phi_j,\ \ \forall \ k\geq
2.\ee Applied to $\{\psi_k\}$ and $\{\varphi_k\}$, the summation in
the right hand side of equation (\ref{4.27}) converges as
$k\rightarrow \infty$, thus, $\tilde d_k$ and $\hat d_k$ decay as
$k^{-3}$. Clearly, $\tilde d_k$ and $\hat d_k$ satisfy the
requirement of Lemma~\ref{l4.2} and for some constants $C_1$, $C_2$,
\be\label{4.28}C_1k^{-3}\leq \tilde d_k,\ \ \hat d_k\leq C_2k^{-3}\
\ \forall\ k\geq 1.\ee

By Lemma~\ref{l4.1}, Lemma~\ref{l4.2} and equation (\ref{4.28}), we
have
\[C_1k^{-3}\leq \tilde d_k=\lim_{t\rightarrow\infty}\frac{\tilde D_k(t)}{t}\leq \liminf_{t\rightarrow\infty}\frac{\overline D_k(t)}{t}
\leq\limsup_{t\rightarrow\infty}\frac{\overline
D_k(t)}{t}\leq\lim_{t\rightarrow\infty}\frac{\hat D_k(t)}{t}=\hat
d_k\leq C_2 k^{-3}\] for all $k\geq 1$. \QED

\section{Proofs of Corollaries~\ref{th3} and \ref{th4}}

\renewcommand{\theequation}{5.\arabic{equation}}
\setcounter{equation}{0}

In this section, we prove Corollaries~\ref{th3} and \ref{th4}.
Because the basic approach is the same as we have used in the proof
of Theorem~\ref{th1}, we only give out a sketch.

For the process $\{G^\a_t:t\geq 1\}$, $0\leq\a<1$, denote by $e_t$
the number of edges in $G^\a_t$ and $a_t=e_{t+1}-e_t$ none the less.

{\it Sketch of the proof of Corollary~\ref{th3}}: For simplicity, we
only deal with the special case of $\mu=\zeta$.

Firstly, it is straightforward to check that Lemmas~\ref{l2.1},
\ref{l2.2} and \ref{l2.3} hold for $e_t$. Then the recurrence of
$\overline D_k(t)$ can be derived as
\be\label{5.1}\left\{\begin{array}{rl}&\hskip-4mm\overline
D_k(t+1)=\overline D_k(t)+\(\Dp\frac{\a(k-1)}{2}+(1-\a)\mu\)\Dp\frac
{\overline D_{k-1}(t)}{t}-\(\Dp\frac{\a k}{2}+(1-\a)\mu\)\Dp\frac
{\overline D_{k}(t)}{t}\\[3mm]&\hskip 18mm+A_k(t)+g^\a_k(t),\hskip 5mm t+1\geq k\geq 0,\ t\geq 1;\\[4mm]
&\hskip-4mm\overline D_0(1)=0;\hskip 3mm\overline D_1(1)=2;\hskip
3mm\overline D_k(t)=0,\ k> t\geq 1;\hskip 3mm\overline D_{-1}(t)=0,\
t\geq 1.\end{array}\right.\ee where $A_k(t)$ satisfying
(\ref{4.6'}), $g^\a_k(t)=0, \ \forall\ t\leq k-1$ and \bena
&&\hskip-4mmg^\a_k(t)=\P(a_t=k)\nonumber\\[3mm]
&&\hskip-4mm=\a\P\(\Dp\sum_{i=0}^t B\(1,\Dp\frac{\mu
d^\a_{x_i}(t)}{2e_t}\)=k\)+(1-\a)\P\(B\(t+1,\Dp\frac{\mu}{t+1}\)=k\)\nonumber\\[3mm]
&&\hskip-4mm=:\a f^\a_k(t)+(1-\a)\bar f_k(t),\ \ \forall\ t\geq
k.\label{5.2}\eena


In the case of $0\leq\a<1$, we have
\be\label{5.3}\liminf_{t\rightarrow\infty}g^\a_k(t)\geq
(1-\a)\lim_{t\rightarrow\infty}\bar
f_k(t)=(1-\a)\Df{\mu^k}{k!}e^{-\mu},\ \ \forall\ k\geq 0,\ee then,
there exists some $\rho>0$ such that (\ref{4.0}) holds. Note that
here we get such $\rho$ from (\ref{5.3}), but in case of $\a=1$, we
get it from the existence of $\lim_{t\rightarrow\infty}\P(a_t=0)$,
which depends on the degree bounds given in Lemma~\ref{l3.1}.

In case of $\a>0$, let $n(\a)=3+\lfloor2/\a\rfloor$, where $\lfloor
2/\a\rfloor$ be the integer part of $2/\a$. It is straightforward to
check that \be\label{5.4}g^\a_0(t)\leq e^{-\mu},\ \forall\
t\geq1\hb{ and }g^\a_k(t)\leq \Df{(\mu\vee 1)^{n(\a)}\times
n(\a)!}{k^{n(\a)}}, \ \forall\ k\geq 1,\ t\geq 1.\ee
Define $\{\psi_k\}$ and $\{\varphi_k\}$ as \[\psi_k=\left\{\ba{rll}&\hskip-4mm0,\hskip3mm &k\geq 1,\\[5mm]
&\hskip-4mm\rho,\hskip3mm &k=0;\ea\right.\hb{ and}\hskip3mm
\varphi_k=\left\{\ba{rll}&\hskip-4mm C(\a)k^{-n(\a)},\hskip3mm
&k\geq 1,\\[5mm]
&\hskip-4mm e^{-\mu},\hskip3mm &k=0,\ea\right.\] with
$C(\a)=(\mu\vee 1)^{n(\a)}\times n(\a)!$. Then define
\[\psi_k(t)=\left\{\ba{rll}&\hskip-4mm\psi_k,\hskip3mm &t\geq \(\Df{\a k}{2}+(1-\a)\mu\)\vee k,\\[5mm]
&\hskip-4mm g^\a_k(t),\hskip3mm &1\leq t<\(\Df{\a k
}{2}+(1-\a)\mu\)\vee k;\ea\right.\] and
\[\varphi_k(t)=\left\{\ba{rll}&\hskip-4mm\varphi_k,\hskip3mm &t\geq \(\Df{\a k}{2}+(1-\a)\mu\)\vee k,\\[5mm]
&\hskip-4mm g^\a_k(t),\hskip3mm &1\leq t<\(\Df{\a k
}{2}+(1-\a)\mu\)\vee k.\ea\right.\] Thus we have
\be\label{5.5}\psi_k(t)\leq g^\a_k(t)\leq\varphi_k(t),\ \ \forall \
k\geq 0,\ t\geq 1.\ee

Let $\tilde D_k(t)$ and $\hat D_k(t)$ be the solutions of the
recurrences obtained from (\ref{5.1}) with $g^\a_k(t)$ substituted
by $\psi_k(t)$ and $\varphi_k(t)$ respectively. Namely
\[\left\{\begin{array}{rl}&\hskip-4mm\tilde
D_k(t+1)=\tilde D_k(t)+\(\Dp\frac{\a(k-1)}{2}+(1-\a)\mu\)\Dp\frac
{\tilde D_{k-1}(t)}{t}-\(\Dp\frac{\a k}{2}+(1-\a)\mu\)\Dp\frac
{\tilde D_{k}(t)}{t}\\[3mm]&\hskip 18mm+A_k(t)+\psi_k(t),\hskip 5mm t+1\geq k\geq 0,\ t\geq 1;\\[4mm]
&\hskip-4mm\tilde D_0(1)=0;\hskip 3mm\tilde D_1(1)=2;\hskip
3mm\tilde D_k(t)=0,\ k> t\geq 1;\hskip 3mm\tilde D_{-1}(t)=0,\ t\geq
1;\end{array}\right.\] and
\[\left\{\begin{array}{rl}&\hskip-4mm\hat
D_k(t+1)=\hat D_k(t)+\(\Dp\frac{\a(k-1)}{2}+(1-\a)\mu\)\Dp\frac
{\hat D_{k-1}(t)}{t}-\(\Dp\frac{\a k}{2}+(1-\a)\mu\)\Dp\frac
{\hat D_{k}(t)}{t}\\[3mm]&\hskip 18mm+A_k(t)+\varphi_k(t),\hskip 5mm t+1\geq k\geq 0,\ t\geq 1;\\[4mm]
&\hskip-4mm\hat D_0(1)=0;\hskip 3mm\hat D_1(1)=2;\hskip 3mm\hat
D_k(t)=0,\ k> t\geq 1;\hskip 3mm\tilde D_{-1}(t)=0,\ t\geq
1.\end{array}\right. \] Then Lemma~\ref{l4.1} holds and we have
\be\label{5.7}\tilde D_k(t)\leq \overline D_k(t)\leq\hat D_k(t),\ \
\forall \ k\geq -1,\ t\geq 1.\ee

Define the two recurrences with respect to $\{\psi_k\}$ and
$\{\varphi_k\}$ respectively as
\[\left\{\ba{rl}&\hskip-4mm\tilde
d_k=\Dp\(\frac{\a(k-1)}{2}+(1-\a)\mu\)\tilde
d_{k-1}-\Dp\(\frac {\a k}2+(1-\a)\mu\)\tilde d_k+\psi_k,\ \ \ k\geq 0,\\[5 mm]
&\hskip-4mm\tilde d_{-1}=0;\ea\right.\] and
\[\left\{\ba{rl}&\hskip-4mm\hat
d_k=\Dp\(\frac{\a(k-1)}{2}+(1-\a)\mu\)\hat
d_{k-1}-\Dp\(\frac {\a k}2+(1-\a)\mu\)\hat d_k+\varphi_k,\ \ \ k\geq 0,\\[5mm]
&\hskip-4mm\hat d_{-1}=0.\ea\right.\]Then Lemma~\ref{l4.2} holds,
namely, under the condition that $\tilde d_k\leq C/k$ (resp. $\hat
d_k\leq C/k$) for some constant $C$ and $k\geq 1$, there exists
constant $M_1$ (resp. $M_2$) such that \be\label{5.8}\left|\tilde
D_k(t)-t\tilde d_k\right|\leq M_1t^{4/5} \(\hb{resp.} \ \left|\hat
D_k(t)-t\hat d_k\right|\leq M_2t^{4/5}\)\ee for all $k\geq -1$ and
$t\geq 1$.

Finally, it suffices to solve the recurrence in $k$ with the form
\be\label{5.9}\left\{\ba{rl}&\hskip-4mm
d_k=\Dp\(\frac{\a(k-1)}{2}+(1-\a)\mu\)
d_{k-1}-\Dp\(\frac {\a k}2+(1-\a)\mu\) d_k+\phi_k,\ \ \ k\geq 0,\\[5 mm]
&\hskip-4mm d_{-1}=0,\ea\right.\ee where $\{\phi_k:k\geq 0\}$ be a
serial of nonnegative numbers. Clearly, recurrence (\ref{5.9}) can
be solved as: $d_{-1}=0$, $d_0=\Dp\frac{2}{b\a}\phi_0$ and
\be\label{5.10}d_k=\prod_{j=1}^k\(1-\frac{\beta}{j+b}\)\(\sum_{i=1}^k\Df{1}{\prod_{j=1}^i(1-\frac{\beta}{j+b})}\frac{2}{(i+b)\a}\phi_i
+\frac 2{b\a}\phi_0\), \hb{ for }\ k\geq 1,\ee where $\beta=1+2/\a$
and $b=2/\a+2(1-\a)\mu/\a$. Applying to $\{\psi_k\}$ and
$\{\varphi_k\}$, the summation term in the right side of equation
(\ref{5.10}) converges as $k\rightarrow\infty$, this implies that
$\tilde d_k$, $\hat d_k$ decay as $k^{-\beta}$. In particular, for
some positive constants $C_1^\a$ and $C^\a_2$, \be\label{5.11}C_1^\a
k^{-\beta}\leq \tilde d_k,\ \hat d_k\leq C_2^\a k^{-\beta},\ \
\forall\ k\geq 1.\ee Corollary~\ref{th3} follows immediately from
(\ref{5.7}), (\ref{5.8}) and (\ref{5.11}).\QED

{\it Sketch of the proof of Corollary~\ref{th4}}: In the case of
$\a=0$, the recurrence of $\overline D_k(t)$ can be derived as
\be\label{5.12}\left\{\begin{array}{rl}&\hskip-4mm\overline
D_k(t+1)=\overline D_k(t)+\zeta\Dp\frac {\overline
D_{k-1}(t)}{t}-\zeta\Dp\frac
{\overline D_{k}(t)}{t}+\bar A_k(t)+g^0_k(t),\\[4mm] & \hskip 63mm t+1\geq k\geq 0,\ t\geq (\zeta-1)\vee 1;\\[4mm]
&\hskip-4mm\overline D_0(1)=0,\hskip 3mm\overline
D_1(1)=2;\hskip3mm\overline D_k(t)=0,\ k> t\geq 1;\hskip
3mm\overline D_{-1}(t)=0,\ t\geq
1;\\[4mm]
&\hskip-4mm\overline D_{k'}(k)=0,\hskip 3mm\overline
D_k(k)=(k+1),\hskip 3mm 0\leq k'<k,\ 1<k<\zeta,\end{array}\right.\ee
where
\[\bar A_k(t)=\Df{\zeta\(\overline D_k(t)-\overline D_{k-1}(t)\)}{t(t+1)}\] and $g^0_k(t)=\bar f_k(t)$, which is given in
(\ref{5.2}) with the parameter $\mu$ replaced by $\zeta$. Note that
the last line in (\ref{5.12}) comes from the fact $G^0_t$ is a
complete graph while $t<\zeta$.

It is clear that $|\bar A_k(t)|\leq {(2\zeta)}/t$ and then satisfies
(\ref{4.6'}), i.e., for some $R_1>0$,
\[|\bar A_k(t)|\leq R_1 t^{-1/5},\ \ \forall \ t\geq 1, \ k\geq 0.\]
For the term $g^0_k(t)=\bar f_k(t)$, we have
\[\lim_{t\rightarrow\infty}\bar f_k(t)=\Df{\zeta^k}{k!}e^{-\zeta},\ \forall\
k\geq 0;\] on the other hand,
\[\bar f_k(t-1)=\binom{t}{k}\(\frac \zeta t\)^k\(1-\frac \zeta
t\)^{t-k}\leq\Df{t(t-1)\cdots(t-k+1)}{(t-\zeta)^k}\Df{\zeta^k}{k!}e^{-\zeta},\]
for all $t\geq \zeta\vee 2$, and $1\leq k\leq t$, this implies that
 \be\label{5.6}\bar f_k(t)\leq C(0)
\Df{\zeta^k}{k!}e^{-\zeta},\ \ \hb{for all } k\geq 1 \hb{ and }t\geq
1\ee for some constant $C(0)>0$.

Now, by (\ref{5.3}), we choose $\rho>0$ satisfying (\ref{4.0}) and
define $\{\psi_k\}$, $\{\varphi_k\}$ as
\[\psi_k=\left\{\ba{rll}&\hskip-4mm0,\hskip3mm &k\geq 1,\\[5mm]
&\hskip-4mm\rho,\hskip3mm &k=0;\ea\right.\hskip3mm
\varphi_k=\left\{\ba{rll}&\hskip-4mm
C(0)\Df{\zeta^k}{k!}e^{-\zeta},\hskip3mm
&k\geq 1,\\[5mm]
&\hskip-4mm e^{-\zeta},\hskip3mm &k=0.\ea\right.\] Then,
Corollary~\ref{th4} follows from the comparing argument used above
and the fact that
\[\sum_{k=0}^\infty\(\frac{1+\zeta}{\zeta}\)^k\phi_k<\infty\] for
$\phi_k=\psi_k$ and $\varphi_k$ respectively.\QED

\br\label{r3} To get the degree distribution by the standard
argument introduced in \cite{CF} and \cite{CFV}, appropriate upper
bounds for $\D_t$, the maximum degree, are always necessary. We
point out that no bounds for $\D_t$ are used in our proofs of
Corollaries~\ref{th3} and \ref{th4}. \er

\section{Application to the Hard Copying Model}

\renewcommand{\theequation}{6.\arabic{equation}}
\setcounter{equation}{0}

It is well known that, besides the BA mechanism, copying is another
mechanism that may lead to power law degree sequence. The basic idea
of copying comes from the fact that a new web page is often made by
copying an old one. A kind of copying models was proposed in Kumar
{\it et al}. \cite{KRRS} to explain the emergence of the degree
power laws in the web graphs. These models are parameterized by a
{\it copy factor} $\a\in (0,1)$ and a constant out-degree $d\geq 1$.
At each time step, one vertex $u$ is added and $d$ out-links are
generated for $u$ as follows. First, an existing vertex $p$ is
chosen uniformly at random; then with probability $1-\a$ the $i^{\rm
th}$ out-link of $p$ is taken to be the $i^{\rm th}$ out-link of
$u$, and with probability $\a$ a vertex is chosen from the existing
vertices uniformly at random to be the destination of the $i^{\rm
th}$ out-link of $u$. It is proved in \cite{KRRS} that the above
copying models possess a power law degree sequence as $d_k\sim
Ck^{-(2-\a)/(1-\a)}$.

In this section, as an application of the comparing argument, we
will introduce a new copying model, here we call it {\it hard
copying} model. Note that another hard copying model is introduced
in \cite{NWC}, which is a mixed model of BA mechanism and hard
copying mechanism.

For fixed $0\leq \a\leq 1$ and $\mu>0$, define random graph process
$\{\bar G^\a_t=(V_t,E_t):t\geq 1\}$ as follows

\vskip 2mm

{\it Time-Step 1.} {Let $\bar G^\a_1$ consists of vertices $x_0,x_1$
and the edge $\langle x_0,x_1\rangle$.}

\vskip 2mm {\it Time-Step $t\geq$ 2.} \begin{enumerate}
  \item with probability $\a$, we generate vertex $ x_{t} $
by copying an existing vertex $x_i$, $0\leq i\leq t-1$ from
$V_{t-1}$ uniformly at random. Note that in this case, all neighbors
of $x_t$ are those of the copied vertex $x_{i}$;
 \item with probability $1-\a$, we add a new vertex
$x_t$ to $\bar G^\a_{t-1}$ and then add random edges incident with
$x_t$ in the {\it classical} manner: for any $0\leq i\leq t-1$, edge
$\langle x_i,x_t\rangle$ is added independently with probability
$(\mu\wedge t) /t$.
\end{enumerate}

As calculated in \cite{NWC}, for the present model, we have
\be\label{6.6}\E(e_{t})=\mu t+O(t^{2\a}).\ee So $e_t=|E_t|$ increase
super-linearly when $\a>1/2$. This makes our model interesting and
deferring from the model introduced in \cite{KRRS}.

Another fact for the present model is that, in any case of $\a$,
\be\label{6.0}\Delta_t\leq t,\ \ \forall \ \ t\geq 1,\ee where
$\Delta_t$ be the maximum degree of $\bar G^\a_t$.

In the case of $\a=0$, $\{\bar G^0_t:t\geq 1\}$ is just
$\{G^0_t:t\geq 1\}$ and its degree sequence is given in
Corollary~\ref{th4} (with $\mu$ in place of $\zeta$). In the case of
$\a=1$, we get a pure hard copying model and, using (\ref{6.0}), the
recurrence of $\overline{D}_k(t)$ can be derived as
\be\label{6.1}\left\{\begin{array}{rl}&\hskip-4mm\overline
D_k(t+1)=\overline D_k(t)+(k-1)\(\Dp\frac {\overline
D_{k-1}(t)}{t+1}-\Dp\frac {\overline D_{k}(t)}{t+1}\),\hskip 3mm t+1\geq k\geq 1;\\[5mm]
&\hskip-7mm\hskip 3mm\overline D_1(1)=2;\hskip 3mm\overline
D_k(t)=0,\ k>t\geq 1;\hskip 3mm\overline D_{0}(t)=0,\ t\geq 1.
\end{array}\right.\ee
By (\ref{6.1}), it is straightforward to prove by induction (in $t$)
that, there exists some $M>0$ such that
\be\label{6.5}\overline{D}_k(t)\leq M(t+1)^{1/2},\ \ \forall\ \
k\geq 1\ee for all $t\geq 1$. Thus we obtain a degenerated degree
distribution as follows.

\bp\label{p1} For any $K\geq 1$, we have \be\label{6.3}
\lim_{t\to\infty}\Df{\sum_{k=1}^K\overline{D}_k(t)}{t+1}=0;\ee
furthermore, for any $\epsilon>0$, we have
 \be\label{6.4}
\lim_{t\to\infty}\Dp\P\(\sum_{k=1}^K{D}_k(t)\geq
\epsilon(t+1)\)=0.\ee\ep

For the case of $0<\a< 1$, using (\ref{6.0}) again, the recurrence
of $\overline{D}_k(t)$ can be derived as
\be\label{6.2}\left\{\begin{array}{rl}&\hskip-4mm\overline
D_k(t+1)=\overline D_k(t)+[\a (k-1)+(1-\a)\mu]\(\Dp\frac {\overline
D_{k-1}(t)}{t+1}-\Dp\frac {\overline D_{k}(t)}{t+1}\)+(1-\a)\bar f_k(t),\\[4mm] & \hskip 63mm t+1\geq k\geq 0,\ t\geq (\mu-1)\vee 1;\\[4mm]

&\hskip-4mm\overline D_k(t+1)=\overline D_k(t)+[\a
(k-1)+(1-\a)(t+1)]\(\Dp\frac {\overline
D_{k-1}(t)}{t+1}-\Dp\frac {\overline D_{k}(t)}{t+1}\)+(1-\a)\bar f_k(t),\\[4mm] & \hskip 63mm t+1\geq k\geq 0,\ 1\leq t<(\mu-1);\\[4mm]&\hskip-4mm\overline D_0(1)=0,\hskip 3mm\overline
D_1(1)=2;\hskip3mm\overline D_k(t)=0,\ k> t\geq 1;\hskip
3mm\overline D_{-1}(t)=0,\ t\geq 1;\end{array}\right.\ee where $\bar
f_k(t)$ is given in (\ref{5.2}).

By the comparing argument developed in Section 4, we can solve
(\ref{6.2}) and obtain the following result.

\bt\label{th5} For any $0<\a<1$ and $\mu>0$, there exists positive
constants $\bar C^\a_1$ and $\bar C^\a_2$ such that
 \be\label{1.3}
\bar C^a_1k^{-1/\a}\leq\liminf_{t\rightarrow\infty }\Df{\overline
{D}_k(t)}{t}\leq \limsup_{t\rightarrow\infty }\Df{\overline
{D}_k(t)}{t}\leq \bar C^\a_2 k^{-1/\a}\ee for all $k\geq 1$.\et

Theorem~\ref{th5} provides an interesting result: in the case of
$e_t$ increasing super-linearly, {\it i.e.} $\a>1/2$, the model
processes power law degree sequence, furthermore, the inverse power
lies in interval (1,2], which is never considered in previous
literature.

\br\label{r4} We note here that, except for (\ref{6.0}), no bounds
for $e_t$ and $\D_t$ are used in our proof of Theorem~\ref{th5}.
Clearly, (\ref{6.0}) holds for all models studied in this paper, and
(\ref{6.0}) implies \be\label{6.8}\overline{D}_k(t)=0,\ \ \forall \
k>t\geq 1.\ee In fact, (\ref{6.8}) is a key evidence to ensure
Lemmas~\ref{l4.1} amd \ref{l4.2} in the comparing argument. \er

\end{document}